\documentclass[11pt]{article}
\def\today{June 6, 2005}
\usepackage{amsmath,amssymb}

\newtheorem{thm}{Theorem}[section]
\newtheorem{co}[thm]{Corollary}
\newtheorem{lemma}[thm]{Lemma}

\newtheorem{assumption}[thm]{Assumption}

\newtheorem{pr}[thm]{Proposition}

\newtheorem{assu1}[thm]{Assumption}

\newtheorem{definition}[thm]{Definition}

\newtheorem{example}[thm]{Example}

\newtheorem{remark}[thm]{Remark}

\newtheorem{protocol}[thm]{Protocol}

\newcommand{\N}{\mathbb{N}}
\newcommand{\R}{\mathbb{R}}

\newcommand{\sgn}{\mathrm{sgn}\,}

\newcommand{\vier}[4]{\left[ \begin{array}{cc}
                   #1 & #2 \\ #3 & #4 \end{array} \right]}
\newcommand{\Section}[1]{\section{#1}\setcounter{equation}{0}}
\newcommand{\Sectionn}[1]{\section*{#1}}

\newcommand{\openbox}{\leavevmode
  \hbox to.77778em{%
    \hfil\vrule
  \vbox to.675em{\hrule width.6em\vfil\hrule}%
  \vrule\hfil}}
\newcommand{\proofname}{Proof}

\newenvironment{proof}[1][\proofname]{\par\normalfont
  \trivlist\item[\hskip\labelsep\itshape #1:]\ignorespaces
  }{\hspace*{1cm}\hspace*{\fill}\openbox \medskip\endtrivlist}

\title{Existence of a Limiting Distribution for the Binary GCD
Algorithm~\thanks{The author was supported in part by SNF grant
200021-103683. }}%
\date{\today}%
\author{G\'erard Maze \\
{\small {\em e-mail:\/} gmaze@math.unizh.ch \vspace{-1mm} }\\ 
{\small Mathematics Institute\vspace{-1mm}}\\ 
{\small University of Z\"urich\vspace{-1mm}}\\ 
{\small Winterthurerstr 190, CH-8057 Z\"urich, Switzerland }
\vspace{3mm} }

\begin{document}\maketitle
\thispagestyle{empty}

\begin{abstract}
  In this article, we prove the existence and uniqueness of a certain
  distribution function on the unit interval. This distribution
  appears in Brent's model of the analysis of the binary gcd
  algorithm. The existence and uniqueness of such a function has been
  conjectured by Richard Brent in his original paper \cite{brent}.
  Donald Knuth also supposes its existence in \cite{knuth} where
  developments of its properties lead to very good estimates in
  relation with the algorithm. We settle here the question of
  existence, giving a basis to these results, and study the
  relationship between this limiting function and the {\it binary
  Euclidean operator} $B_2$, proving rigorously that its derivative
  is a fixed point of $B_2$.\\

  \noindent{\bf Subject Class:} Primary 68W40, 47H10. Secondary 68Q25
\end{abstract}

\Section{Introduction}

If $u$ and $v$ are positive integers, their {\it greatest common
divisor} (gcd), written $\gcd(u,v)$ in the sequel, is the largest
integer that divides them both. This integer can be computed
efficiently using a method discovered more than 2200 years ago:
Euclid's algorithm. Quoting D.Knuth \cite{knuth}, this algorithm is
the ``grand-daddy'' of all algorithms, because it is the oldest
nontrivial algorithm that has survived to the present day. It is
however not the best way to find greatest common divisors when dealing
with modern computers.  In fact, another algorithm, the so-called {\it
binary gcd algorithm}, created by J.Stein \cite{stein}, requires no
division but only subtractions, parity testings and halving of even
numbers (which correspond to shifts in binary notations). These
procedures are essentially free when compared to the computational
cost of divisions.

The idea of the binary gcd algorithm is basically the following: given
two positive integers $u$ and $v$, if halving both numbers is possible
at most $k$ time, do it, keeping the values of $u$ and $v$
updated. Then repeat the following procedure until both number are
equal, say to $l$: subtract the smaller from the greater and when the
result is even, divide it by the largest power of 2 possible. The gcd
of $u$ and $v$ is then $l \cdot 2^k$. This repeated loop will be
referred as a ``subtract-and-shift cycle'' in the sequel.

The behavior of the binary gcd algorithm is interesting in
several ways. On the one hand, it is always important to know the
worst case and average case of an algorithm, just from a
practical point of view, and this is even more important when the
considered algorithm has such a wide application. On the second
hand, the machinery elaborated in order to understand the average
behavior of the algorithm has led to a deep understanding of it,
giving answers as well as rising new questions.

In our case, the worst case the binary gcd algorithm may have to
face is a total number of subtractions equal to $1+\lfloor \log_2
\max(u,v) \rfloor)$, see, e.g., \cite{knuth}.

The exact determination of the average behavior of the binary gcd algorithm is
however much more complex than the analysis of its worst case scenario. Two
models have been proposed in order to study and analyze the expected behavior
of the algorithm. We first describe the model created by R.Brent and gives a
short description of the model created by B.Vall\'ee at the end of this
introduction.

The first model was created by Richard Brent in 1976 \cite{brent}. In his
work, R.Brent exhibits a dynamical system describing the binary Euclidean
algorithm and provides an heuristic proof of the analysis of the algorithm.
This dynamical system is described by the binary Euclidean operator $B_2$, see
(\ref{B}) below, that transforms the density associated to the algorithm,
step-by-step.  However, the operator $B_2$ is difficult to analyze, and the
question of convergence was left as a conjecture. This approach also suffers
from the fact that it lies on an unproven connection between a discrete and a
continuous model, see \cite{brent} for more details concerning this last
point and \cite{brent2} for a description of the situation 25 years later.

We now describe this model. Suppose that both $u$ and
$v$, with $u>v$, are odd, which is the case after each subtract-and-shift
cycle.  Every subtract-and-shift cycle forms $u-v$ and shifts this quantity
right until obtaining an odd number $u'$ that replace $u$. Under random
conditions, one would expect to have $u'=(u-v)/2^m$ with probability $2^{-m}$.
This is the heart of Brent's hypothesis. In his model, we suppose that $u$ and
$v$ are essentially random, except that they are odd and their ratio $v/u$ has
a certain probability distribution.  Let $g_n$ be the probability that
$\min(u,v)/\max(u,v)$ is greater or equal to $x$ after $n$
subtraction-and-shift cycles have been performed under this assumption. Then
the sequence of functions $\{g_n\}_{n \in \N}$ satisfies the following
recurrence relation \cite{brent,knuth}:
$$
g_0(x)  =  1-x\; , \; \; g_{n+1}(x) =  F(g_n)(x)
$$
where, for all $h \in C([0,1])$,
\begin{equation}\label{F}
F(h)(x) = \sum_{k \geq 1} 2^{-k} \left( h \left(\frac{x}{x+2^k}
\right) - h \left(\frac{1}{1+2^kx} \right) \right), \; x \in
[0,1].
\end{equation}

In the sequel, we will denote by $F_n(h)$ the partial sums of the
above series. Note that the operator $F$ is linear, bounded, since
$||F||_{\infty} \leq 2$, and that the series converges uniformly for
any $h$ in $C([0,1])$. Computational experiments led Richard Brent to
conjecture that the functions $g_n$ converge uniformly to a limiting
distribution $g_{\infty}$. Under this conjecture, the function
$g_{\infty}$ satisfies the equality
\begin{equation*}
 g_{\infty}(x) = \sum_{k \geq 1} 2^{-k} \left(
g_{\infty} \left(\frac{x}{x+2^k} \right) - g_{\infty}
\left(\frac{1}{1+2^kx} \right) \right)
\end{equation*}
and provides the following estimate. If
$$
b = 2+ \int_{0}^{1}\frac{g_{\infty}(x)}{(1-x)\ln 2} dt =
2.83297657... ,
$$
then the expected number of subtract-and-shift cycles in the binary gcd
algorithm with starting values $u$ and $v$ is $\ln(uv)/b$.

In this article, we settle this conjecture by proving that the
sequence $\{g_n\}_{n \in \N}$ converges uniformly towards a function
$g_{\infty}$. In order to do so, we first prove that every element of
the sequence is convex and differentiable over $]0,1]$. Then we
exhibit a compact set of the Banach space $(C([0,1]),||.||_{\infty})$
which contains $g_1$ and which is left invariant by the operator $F$
defined by (\ref{F}). This fact assures the existence of accumulation
points of the sequence $\{g_n\}_{n \in \N}$, and therefore proves the
existence of at least one fixed point of the operator $F$. We prove
the uniqueness via an argument based on the sequence of
derivatives $\{g_n'\}_{n \in \N}$. On the way, we study the behavior
of this sequence with respect to the binary Euclidean operator $B_2$,
proving that the sequence converges to the unique fixed point of $B_2$
in the $L^1$-norm.

The present work provides a proof that the dynamical system studied by
Brent possesses indeed a unique limiting distribution. However, it
does not shed a new light on the validation of the continuous model.
In other words, it makes legitimate the work of R.Brent on the
analysis of the binary gcd algorithm, without validating his model. It
also answers a 47-points question of D.Knuth \cite[p.355, question
32]{knuth}, who grades the problems of \cite{knuth} on a
``logarithmic'' scale from 0 to 50.

The second work we were referring to is due to Brigitte Vall\'ee \cite{vallee1}
who brings another look to the situation and leads to a successful analysis
using rigorous ``dynamical'' methods. These methods are also the basis for the
analysis of several others algorithms \cite{vallee1,vallee2}. In her work,
B.Vall\'ee studies the operator $V_2$ which describes a slightly different
dynamical system. The operator $V_2$ transforms the density associated to the
algorithm where all the subtract-and-shift cycles are gathered together as
long as the sign of $u-v$ is constant. As a consequence, the operator $V_2$ is
easier to analyze. B.Vall\'ee shows that the operator $V_2$ possesses a unique
fixed point in some Hardy space and presents a spectral gap. She also proves,
based on this spectral gap and with the help of a Tauberian theorem, the
connection between the discrete and the continuous model. Quoting D.Knuth,
``her methods are sufficiently different that they are not yet known to
predict the same behavior as Brent's heuristic model. Thus the problem of
analyzing the binary gcd algorithm [...]  continues to lead to ever more
tantalizing questions of higher mathematics''.

Not surprisingly, there is a connection between the two operators $B_2$ and
$V_2$. We refer the interested reader to \cite{brent2} for further details
regarding this connection.
\\

We will use the notation $||.||_{\infty}$ and $||.||_1$ for the
supremum norm and the $L^1$-norm of functions defined over $[0,1]$,
and $\log_2$ for the logarithm in base $2$. Let us recall that a
series of function $\sum_{n>0} h_n(x)$ verifies the so-called
Weierstrass criterion (see \cite[III.4]{wanner}) over a subset $A$ of
$\R$ if we have
\begin{equation}\label{abs}
\sum_{n>0} \;\; \sup_{x \in A}|h_n(x)| < \infty.
\end{equation}

Let us mention that the notation for our $g_n$ and $g_{\infty}$ are
different in both \cite{brent} and \cite{knuth}. Brent \cite{brent}
uses $F_n$ and $F_{\infty}$ and Knuth \cite{knuth} uses $G_n$ and $G$.

\Section{Convexity and Regularity}

We prove in this section that the elements of the sequence
$\{g_n\}_{n \in \N}$ are convex and decreasing functions over $[0,1]$ and
differentiable over $]0,1]$.\\

Let $m$ be a $2\times 2$ matrix with real coefficients $a,b,c$
and $d$ . Such a matrix acts naturally on $\R$ via
$$
m(x) = \vier abcd (x) = \frac{ax+b}{cx+d}
$$
and this action satisfies $(m_1 \cdot m_2)(x) = m_1(m_2(x))$ for
all pairs $m_1,m_2$ of $2 \times 2$ matrices, where $\cdot$ is the
usual matrix product. From now on, we will identify a matrix $m$ with
the real function associated to it. Let us define the set $M$ as
follows:
$$ 
M=\left\{m:[0,1] \rightarrow \R \, \left| \right. \, m(x) =
\frac{ax+b}{cx+d}\, \mbox{ with } \, a,b \geq 0, \; c,d > 0 \, , \,
ad-bc \neq 0 \right\}.
$$
Note that for any element $m$ of $M$, $\sgn (m):=
\mbox{signum}(ad-bc)$ is well-defined, since a common factor at the
denominator and the numerator of $m$ does not affect the sign of
$ad-bc$. This function satisfies the equality $\sgn (m_1 \circ m_2) =
\sgn (m_1) \cdot \sgn (m_2)$ for all pairs $m_1,m_2$ in $M$.  Let us
also define the set $S$ as the set of all series $\sum_{i \in \N}
\varepsilon_i m_i$, where $m_i \in M$, satisfying the following three
points:
\begin{enumerate}
\item $\varepsilon_i = \pm 1$ and $\varepsilon_i \cdot \sgn(m_i) <0 \,
  , \, \forall i \in \N$.  \item The series verifies the 
  Weierstrass criterion over $[0,1]$, i.e.,
  \begin{equation}\label{abs3}
    \sum_{i \in \N} \varepsilon_i m_i \in S \;\; \Longrightarrow \;\;
     \,\, \sum_{i \in \N} ||m_i||_{\infty} < + \infty.
  \end{equation}
\item The following series converges:
  \begin{equation}\label{inv}
    \sum_{i \in \N} \frac{|a_id_i-b_ic_i|}{c_id_i} < +\infty\, , \,
    \mbox{ where } 
    m_i(x) = \frac{a_ix+b_i}{c_ix+d_i}. 
  \end{equation}
\end{enumerate}

Note that the series (\ref{inv}) is well-defined since a common
factor at the denominator and the numerator of $m_i$ does not affect
the terms of the series. A typical element $g$ of $S$ can be written as
\begin{equation}\label{g}
g(x)=\sum_{i\in \N} \frac{a_ix+b_i}{c_ix+d_i}.
\end{equation}

For sake of clarity, let us recall two facts about series of
functions: first, if a series of functions satisfies the Weierstrass
criterion (\ref{abs}) on a set, then it converges uniformly and
absolutely on it, and the limit does not depend on any permutation of
the sum. This result applies also for double sums. Second, if the
derivatives of the partial sums of a convergent series of function
converges uniformly then the series is differentiable and its
derivative is the limit of the derivatives of the partial sums, see,
e.g., Thm. 2.13, Thm. 2.9, Thm. 4.3 and Thm. 6.18 of \cite{wanner}.

The definition of the set $S$ takes its roots in the following two
lemmas, which are the keystones of the article.

\begin{lemma}\label{conv}
Every function $g$ in $S$ is a convex, decreasing, continuous function
over $[0,1]$ and continuously differentiable over any compact set of
$]0,1]$.
\end{lemma}

\begin{proof}
A function $m$ in $M$ is convex and decreasing if and only if $\sgn(m)
< 0$. Indeed, we have
$$
\left(\frac{ax+b}{cx+d}\right)'=  \frac{ad-bc}{(cx+d)^2}<0
\;\; \Longleftrightarrow \;\; ad-bc<0,
$$
and
$$
\left(\frac{ax+b}{cx+d}\right)''=-2c \cdot \frac{ad-bc}{(cx+d)^3}>0
\;\; \Longleftrightarrow \;\; ad-bc<0.
$$
Using the first two points of the above definition, any element $g$
in $S$ is a uniform limit of convex, decreasing and continuous
functions, and is therefore convex, decreasing and continuous.  Let us
prove now that any element of $S$ is continuously differentiable over
any compact interval of $]0,1]$. Let $0< \varepsilon < 1$. For $g$ as
in (\ref{g}), the partial sums of $g'$ satisfy
$$
\left(\sum_{i=0}^{N} \frac{a_ix+b_i}{c_ix+d_i}\right)'
 =  \sum_{i=0}^{N} \frac{a_id_i-b_ic_i}{(c_ix+d_i)^2} 
$$
and the definition of $S$ shows that this series satisfies the 
Weierstrass criterion over $[\varepsilon,1]$ since for all $x \in
[\varepsilon,1]$,
$$ 
\frac{|a_id_i-b_ic_i|}{(c_ix+d_i)^2} \leq
\frac{|a_id_i-b_ic_i|}{(c_i \varepsilon+d_i)^2}\leq \varepsilon^{-2}
\frac{|a_id_i-b_ic_i|}{(c_i +d_i)^2}
$$
and  
$$
\frac{1}{(c +d)^2} \leq \frac{1}{4c d} \; \; \; \;  \forall c,d >0,
$$
yields 
$$ 
\varepsilon^{-2}\sum_{i\geq0} \frac{|a_id_i-b_ic_i|}{(c_i+d_i)^2} \leq
\frac{\varepsilon^{-2}}{4} \cdot \sum_{i\geq 0}
 \frac{|a_id_i-b_ic_i|}{c_id_i} < + \infty.
$$ 
Thus, the partial sums of derivative converge uniformly over
$[\varepsilon,1]$ to a limiting function which is the derivative of
$g$. This finishes the proof.
\end{proof}

\begin{lemma}\label{stable}
Let $F_n$ be the partial sums of the series (\ref{F}). If $g:[0,1]
\longrightarrow \R$ is a function in $S$, then $F_n(g) \in S$ for all
$n \in \N$ and $F(g) \in S$.
\end{lemma}

\begin{proof}
Let us define the following particular elements of $M$:
$$
\mu_k (x)=\left[ \begin{array}{ccc} 1 &\;& 0 \\ 1 &\;& 2^k
\end{array} \right](x) = \frac{x}{x+2^k}
\;\; \mbox{ and } \;\; \nu_k(x) =\left[ \begin{array}{ccc} 0 &\;& 1 \\
2^k &\;&1
\end{array} \right](x) = \frac{1}{2^kx+1}.
$$
Note that these functions map the interval $[0,1]$ in itself and
satisfy
$$
\sgn (\mu_k) >0 \, \mbox{ and } \, \sgn (\nu_k)<0.
$$ 
Let us prove that if $g(x)=\sum_{i \in \N} \varepsilon_i m_i(x)$ is
in $S$, then $F(g)$ lies inside $S$. The proof for the partial sums
$F_n(g)$ is similar, although infinite sums might become finite. We have
\begin{eqnarray}
F(g)(x) & = & \sum_{k \geq 1} 2^{-k} \left(g(\mu_k(x))-g(\nu_k(x))
\right) \nonumber\\
& = & \sum_{k \geq 1} 2^{-k} \left(\sum_{i \in \N} \varepsilon_i
m_i(\mu_k(x))-\sum_{i \in \N} \varepsilon_i m_i(\nu_k(x))
\right)\label{double}.
\end{eqnarray}
Based on the Weierstrass criterion (\ref{abs3}), we have
\begin{eqnarray}
||F(g)||_{\infty} & \leq & \sum_{k \geq 1} 2^{-k} \left(\sum_{i \in \N} 
||m_i \circ \mu_k||_{\infty}+\sum_{i \in \N} ||m_i \circ \nu_k||_{\infty} 
  \right)\nonumber \\
& \leq & \sum_{k \geq 1} 2^{-k} \left(\sum_{i \in \N} ||m_i||_{\infty} 
\right.\left.  +\sum_{i \in \N}||m_i||_{\infty} \right)< + \infty 
\label{abs2}
\end{eqnarray}
and therefore the double sums in (\ref{double}) can be rearranged
in any simple sum
$$
F(g)(x)=\sum_{i \in \N} \varepsilon_i M_i(x)
$$
where $\varepsilon_i M_i(x)$ is either of the type $\varepsilon_j
\cdot 2^{-k}m_j(\mu_k(x))$ or of the type $-\varepsilon_j\cdot
2^{-k}m_j(\nu_k(x))$. Clearly, $F(g)$ has the correct structure to be
an element of $S$. We must now prove that this function fulfills the
three points of the definition of the set $S$.  Inequality
(\ref{abs2}) shows that the latter series fulfills the Weierstrass
criterion, directly proving the second point. Since
$$
\varepsilon_j\cdot \sgn(2^{-k} \cdot  (m_j \circ \mu_k)) = \varepsilon_j
\cdot \sgn( m_j) \cdot \sgn(\mu_k) = \varepsilon_j
\cdot \sgn(m_j) <0,
$$
and
$$
 \varepsilon_j\cdot  \sgn( -2^{-k} \cdot (m_j \circ \nu_k)) =
-\varepsilon_j  \cdot \sgn( m_j) \cdot \sgn(\nu_k) =
\varepsilon_j \cdot \sgn(m_j) <0,
$$ 
the first point is verified. Let us check the validity of the third
point. A straightforward computation shows that if $g$ is as in
(\ref{g})
then the analogue series as (\ref{inv}) for $F(g)$ is the following
double series:
$$ 
\sum_{k \geq 1} \sum_{i \in \N}
\frac{2|a_id_i-b_ic_i|}{(c_i+d_i)d_i\cdot2^k}.
$$ 
This series is convergent since the inequality $(c_i+d_i)d_i >
c_id_i $ yields the following estimate:
\begin{eqnarray*}
\sum_{k \geq 1} \sum_{i \in \N}
 \frac{2|a_id_i-b_ic_i|}{(c_i+d_i)d_i\cdot2^k} & < & 
 2 \sum_{k \geq 1} \frac{1}{2^k} \cdot \sum_{i \in \N}
 \frac{|a_id_i-b_ic_i|}{c_id_i}\\
& = & 2 \sum_{i \in \N}
 \frac{|a_id_i-b_ic_i|}{c_id_i} < +\infty.
\end{eqnarray*}
This shows that the third point is fulfilled and the lemma is then
proven.
\end{proof}

\begin{pr}\label{kl}
  Every element of the sequence $\{g_n\}_{n \in \N \setminus \{0\}}$ is
  in $S$. Thus every element of the sequence $\{g_n\}_{n \in \N}$ is a
  convex, continuous and decreasing function over $[0,1]$,
  continuously differentiable over $]0,1]$.
\end{pr}

\begin{proof}
Since $g_0(x)=1-x$, $g_0$ fulfills the conditions of the claim. The
function $g_1$ is as follows
$$ 
g_1(x) = \sum_{k \geq 1} 2^{-k} \left( \frac{1}{1+2^kx} -
\frac{x}{x+2^k} \right) = \sum_{k \geq 1} \left( \frac{1}{2^k+2^{2k}x}
- \frac{x}{2^kx+2^{2k }} \right),
$$ 
and a straightforward computation shows that the function $g_1$
above is an element of $S$. Lemma \ref{stable} shows by induction
that $g_n$ is an element of $S$ and is therefore convex, decreasing,
continuous over $[0,1]$ and continuously differentiable over any
compact subset of $]0,1]$ by Lemma \ref{conv}.
\end{proof}

\Section{Existence of an accumulation point}

In the current section, we prove that the sequence $\{g_n\}_{n
\in \N}$ possesses at least one accumulation point in the Banach
space of continuous function defined over $[0,1]$, with the
supremum norm. Let us define the following two subsets of this
Banach space:
\begin{eqnarray}
K_1 & = & \overline{S} \cap \{ g  |  \, g(0)=1, g(1)= 0 \}, \label{K1}\\
K_2 & = & \{ g \in C([0,1]) \, | \, 1+3/2 \cdot x \log_2 x -5x \leq
g(x) \leq 1-x \}, \label{K2}
\end{eqnarray}
where $\overline{S}$ is the closure of $S$ in the supremum norm and
$\log_2$ is the logarithm in base 2. Note that any element of
$\overline{S}$ is a decreasing, convex and continuous function, being
a uniform limit of such functions. The definition of $K_2$ seems odd
at first sight. The key point is that a function in $K_2$ cannot come
close to 1 with a too steep slope when $x$ goes to 0. We start with
the following proposition:
\begin{pr}
The operator $F$ verifies the following properties:
\begin{enumerate}
\item $F(K_1) \subset K_1$,
\item $F(K_1 \cap K_2) \subset K_2$,
\end{enumerate}
and therefore $F(K_1 \cap K_2) \subset K_1 \cap K_2$.
\end{pr}

\begin{proof}
  The map $F$ being continuous, we have $F(\overline{S}) \subset
  \overline{F(S)}$. Lemma \ref{stable} implies that $ \overline{F(S)}
  \subset \overline{S} $, therefore $F(\overline{S}) \subset
  \overline{S}$. The fact that $F(g)(0)=1$ and $F(g)(1)=0$ when
  $g(0)=1$ and $g(1)=0$ is straightforward. This proves the
  first point.
  
  Suppose $g$ is a function in $K_1 \cap K_2$. The inequality $F(g)(x)
  \leq 1-x$ is obvious since, $F(g)$ being an element of
  $\overline{S}$, is convex and lies below the secant joining
  $(0,1)$ to $(1,0)$. It remains to show that
  $$
  F(g)(x) \geq 1+3/2 \cdot x \log_2 x -5x.
  $$
  Based on the definitions of $F$ and $K_2$, we have
\begin{eqnarray*}
F(g)(x) & = & \sum_{k \geq 1} 2^{-k} \left(
g(x/(x+2^k))-g(1/(1+2^kx))
 \right)\\
& \geq & \sum_{k \geq 1} 2^{-k} \left( 1+\frac{3}{2} \cdot \left(
\frac{x}{x+2^k} \cdot \left(\log_2 x - \log_2 (x+2^k)
\right)\right) \right.\\
& & \left. -5\cdot\frac{x}{x+2^k}
-1+\frac{1}{1+2^kx} \right)\\
& = & \frac{3}{2}\cdot x\log_2 x \cdot  \left( \sum_{k \geq 1} \frac{1}{2^k}
\cdot \frac{1}{x+2^k} \right) - \frac{3}{2}\cdot x \cdot  \left( \sum_{k \geq 1}
\frac{1}{2^k} \cdot \frac{\log_2 (x+2^k)}{x+2^k} \right) \\
& &  -5x \cdot  \left( \sum_{k \geq 1} \frac{1}{2^k} \cdot \frac{1}{x+2^k}
 \right) +  \left(\sum_{k \geq 1} \frac{1}{2^k} \cdot \frac{1}{1+2^kx} \right).
\end{eqnarray*}
Note that, for $x \in [0,1]$, we have
$$
\sum_{k \geq 1} \frac{1}{2^k} \cdot \frac{1}{x+2^k} \leq \sum_{k
\geq 1} \frac{1}{4^k} = \frac{1}{3},
$$
and
$$
\sum_{k \geq 1} \frac{1}{2^k} \cdot \frac{\log_2 (x+2^k)}{x+2^k} \leq
\sum_{k \geq 1} \frac{1}{2^k} = 1.
$$
Based on Mellin's transform, the equality
\begin{equation}\label{P}
\sum_{k \geq 1} \frac{1}{2^k} \cdot \frac{1}{1+2^kx}= 1+x \log_2 x +
x\cdot P(\log_2 x) +\frac{x}{2}- \sum_{k\geq2} (-1)^k
\frac{2^{k-1}}{2^{k-1}-1} x^k
\end{equation}
where
$$
P(y)= \frac{2\pi}{\ln 2} \cdot \sum_{k\geq1} \frac{\sin 2 \pi k
  y}{\sinh (2k\pi^2/\ln 2)}
$$
can be proven. A proof can also be found in \cite[p. 644]{knuth}, where it
appears as the main step in the computation of the function $g_1$.  As
a matter of fact, the function $P(y)$ is small, and can be bounded in
absolute value by $8 \cdot 10^{-12}$, c.f. \cite{knuth}. We will
however only need a far less accurate bound. Since $\sinh(t)>e^t/4$
for $t >\ln 2/2$, we have
$$
|P(y)|<\frac{2\pi}{\ln 2} \cdot \sum_{k\geq1} \left(4e^{-2\pi^2/\ln
    2}\right)^k =\frac{2\pi}{\ln 2} \cdot \frac{4}{e^{2\pi^2/\ln
    2}-4}= 1.5549..\cdot 10^{-11} < 1/4.
$$
For $x \in [0,1]$, the terms of the alternating sums on the
right-hand-side of (\ref{P}) decrease in absolute value. This sum
can therefore be bounded above by its first term, $2x^2$. Since
$x\log_2 x \leq 0$ over $[0,1]$,  the previous estimation of $F(g)$
becomes
\begin{eqnarray*}
F(g)(x) & \geq & \frac{3}{2}\cdot x\log_2 x \cdot \frac{1}{3} -
\frac{3}{2}\cdot x \cdot 1\\
& &  -5 \cdot x \cdot \frac{1}{3} + 1+x \log_2 x - x\cdot \frac{1}{4}
+\frac{x}{2}- 2x^2\\
& = & 1+\frac{3}{2}\cdot x \log_2 x - \frac{35}{12}x -2x^2\\ 
& = &
1+\frac{3}{2}\cdot x \log_2 x - 5x +\underbrace{ \left(\frac{25}{12}x
-2x^2\right)}_{\geq 0}\\ & \geq & 1+\frac{3}{2}\cdot x \log_2 x - 5x.
\end{eqnarray*}
This last estimate finishes the proof of the proposition.
\end{proof}

Note that the proof of the previous proposition also shows that if a
function $g$ is convex and in $K_2$, then $F(g)$ is in $K_2$ as
well. Indeed, the only property needed from $K_1$ in the proof that
$F(K_1 \cap K_2) \subset K_2$ is the convexity of elements in
$K_1$. Let us state this result as a corollary: 

\begin{co}\label{co}
If a function $g:[0,1] \rightarrow \R$ is convex and in $K_2$, then
$F(g)$ is in $K_2$.
\end{co}

We turn now to a result of compactness.
\begin{pr}
  The set $K_1 \cap K_2$ is compact in the Banach space $(C([0,1]),
  ||.||_{\infty})$.
\end{pr}

\begin{proof}
  The set $K_1$ is closed in $(C([0,1]), ||.||_{\infty})$ being the
  intersection of two closed sets. The set $K_2$ is clearly closed as
  well.  Consider the set of H\"older functions over $[0,1]$ (with
  parameter $1/2$). These are the functions $f$ for which
  $$
  N_{1/2}(f)=\sup_{x \neq y} \frac{|f(x)-f(y)|}{|x-y|^{1/2}} <
  \infty .
  $$
  Based on an argument of equicontinuity, it can be verified that
  the set
  $$
  K_{A,B}=\{f \in C([0,1]) \, | \, ||f||_{\infty} \leq A , \,
  N_{1/2}(f)\leq B \}
  $$ 
  is compact in $C([0,1])$ for any $A,B>0$, see e.g.,
  Chap.4, Sect.6 of \cite{folland}.  Let us show that $K_1 \cap K_2
  \subset K_{1,5}$.  The only non-trivial point to be checked is the
  fact that if $g \in K_1 \cap K_2$, then $N_{1/2}(g) \leq 5$. The
  function $g$ being decreasing and convex, we have, for $0<y<x<1$,
  $$
  \frac{|g(x)-g(y)|}{|x-y|^{1/2}} \leq
  \frac{g(0)-g(x-y)}{\sqrt{x-y}}=\frac{1-g(h)}{\sqrt{h}}, \;\; \mbox{
    with } h = x-y>0.
  $$
  Using a property of the elements of $K_2$, we also have
  $$
  \frac{1-g(h)}{\sqrt{h}} \leq \frac{-3/2h \log_2 h +5h}{\sqrt{h}}=
  5\sqrt{h} -3/2 \sqrt{h} \log_2 h.
  $$
  The maximum value of the latter function, defined over $[0,1]$,
  is reached for $h=1$ and therefore
  $$
  \frac{|g(x)-g(y)|}{|x-y|^{1/2}} \leq \left[ 5\sqrt{h} -3/2
    \sqrt{h} \log_2 h \right]_{h=1} = 5.
  $$ 
  Taking the supremum, we obtain the expected result. The set $K_1
  \cap K_2$ being a closed subset of a compact metric space, it is
  itself compact. This proves the proposition.
\end{proof}

The previous two propositions give directly the next corollary, since
any compact set in a metric space satisfies the Bolzano-Weierstrass
condition:

\begin{co}
The sequence $\{g_n\}_{n \in \N}$ possesses at least one accumulation
point in $K_1 \cap K_2$.
\end{co}

\begin{proof}
The function $g_0$ is convex and in $K_2$. By Corollary \ref{co},
$F(g_0)=g_1$ is therefore in $K_2$. This function is also in $K_1$
(see the proof of Proposition \ref{kl}) and thus any element of the
sequence $\{g_n\}_{n \in \N}$ but $g_0$  is in the compact $K_1
\cap K_2$. The conclusion follows by the Bolzano-Weierstrass property.
\end{proof}

\Section{Behavior of the derivatives and uniqueness of the
accumulation point }

In the current section, we prove that the sequence $\{g_n\}_{n \in
  \N}$ in fact possesses only one accumulation point $g_{\infty}$.
This proves that the sequence converges to this well-defined function
in $K_1 \cap K_2$ since a sequence in a compact metric space with only
one accumulation point converges to this point. In order to achieve
this goal, we study the sequence of derivatives $\{g'_n\}_{n \in \N}$
in the topology of the $L^1$-norm over $]0,1]$.  Then, based on a
property of the binary Euclidean operator $B_2$, defined below by
(\ref{B}), we show the uniqueness of the accumulation points of both
the sequences $\{g_n\}_{n \in \N}$ and $\{g'_n\}_{n \in \N}$.


The sequence of derivatives $\{g'_{n}\}_{n \in \N}$ does not converges
uniformly, or at least we do not know it, and therefore nothing tells
us that an accumulation point $g_{\infty}$ possesses a derivative
which is the uniform limit of the derivatives of the
subsequence. However, $g_{n}$ being convex for all $n$, we will see in
the sequel that the sequence of derivatives does converge but in a
weaker topology, the topology of $L^1([0,1])$. Recall that any convex
function $h$ over $[0,1]$ is differentiable almost everywhere in
$[0,1]$ and that it is absolutely continuous \cite{folland}, i.e.,
\begin{equation*}
h(x) = h(0) + \int_0^x h'(t)dt, \;\; \forall x \in [0,1].
\end{equation*}

The following lemma sheds light on the convergence of the
derivatives of convex functions:
\begin{lemma}\label{lem}
Let $\{f_n\}_{n \in \N}$ be a sequence of convex functions defined
over $[a,b]$, differentiable over $]a,b[$, and converging
uniformly to a function $f$. Then:
\begin{enumerate}
\item
If $E$ is the subset of point of $]a,b[$ where $f$ is
differentiable, then for all $x_0$ in $E$, the sequence
$\{f'_n(x_0)\}_{n \in \N}$ converges to $f'(x_0)$.
\item
If the functions $f_n$ are monotone, then the sequence
$\{f'_n\}_{n \in \N}$ converges to $f'$ in the $L^1$-norm.
\end{enumerate}

\end{lemma}
\begin{proof}
Let $x_0 \in E$. There exists $h_0>0$ such that $[x_0-h_0,x_0+h_0]
\subset ]a,b[$. The functions $f_n$ being convex, we have $\forall n
\in \N\, , \; \forall \, 0<h<h_0$,
$$
\frac{f_n(x_0-h)-f_n(x_0)}{h} \leq f'_n(x_0) \leq
\frac{f_n(x_0+h)-f_n(x_0)}{h}.
$$
When $n$ goes to infinity, this leads to the following inequalities
$\forall \, 0<h<h_0$,
$$
\frac{f(x_0-h)-f(x_0)}{h} \leq \liminf_{n} f'_n(x_0) \leq
\limsup_{n} f'_n(x_0) \leq \frac{f(x_0+h)-f(x_0)}{h}.
$$
Taking the limit when $h$ goes to $0$, we finally have
$$
f'(x_0) \leq \liminf_{n} f'_n(x_0) \leq \limsup_{n} f'_n(x_0) \leq
f'(x_0), \mbox{ i.e., } \lim_{n \rightarrow \infty} f'_n(x_0) =
f'(x_0).
$$
In order to prove the second point, note that since $f$ is convex,
the set $[a,b] \setminus E$ has measure $0$, and thus $f_n'$ converges
to $f'$ almost everywhere. Without loss of generality, suppose the
functions $f_n$ are increasing, i.e., $f_n' \geq 0$ almost everywhere.
The functions $f_n$ and $f$ are absolutely continuous, and
thus
$$ 
\lim_{n \rightarrow \infty} \int_a^b |f_n'| = \lim_{n \rightarrow
\infty} \int_a^b f_n' = \lim_{n \rightarrow \infty} \left(
f_n(a)-f_n(b) \right)= f(a)-f(b) = \int_a^b f' = \int_a^b |f'|.
$$ 
A direct application of the dominated convergence theorem shows
that $f'_n$ converges to $f'$ in $L^1$, see also Ex.21, p.57
of \cite{folland}.
\end{proof}

Consider the following linear operator, obtained by taking the
formal derivative of the series  (\ref{F}):
\begin{equation}\label{B}
B_2(h)(x) = \sum_{k \geq 1} \left(\frac{1}{x+2^k} \right)^2 h
\left(\frac{x}{x+2^k} \right) + \left(\frac{1}{1+2^kx} \right)^2 h
\left(\frac{1}{1+2^kx} \right).
\end{equation}

This operator is referred as the ``binary Euclidean operator'' in the
literature. It was first studied by Richard Brent \cite{brent}.
It is not clear at first sight for what class of function the
operator $B_2$ should be defined. If we consider its action on
$L^1([0,1])$, then the operator $B_2$ is a contraction with respect
to the $L_1$-norm:
\begin{eqnarray}
\int_0^1 |B_2(h)(t)|dt & \leq & \sum_{k \geq 1} \int_0^1
\left(\frac{1}{t+2^k} \right)^2
\left| h \left(\frac{t}{t+2^k} \right) \right| dt \nonumber \\
& & + \int_0^1 \left(\frac{1}{1+2^kt} \right)^2 \left|h
\left(\frac{1}{1+2^kt} \right)\right| dt \nonumber \\
& = & \sum_{k \geq 1} 2^{-k} \left(\int_0^{1/(1+2^k)}  |h(y)| dy  +
\int_{1/(1+2^k)}^1 |h(y)| dy \right) \nonumber \\
& = & \int_0^1 |h(y)|dy. \label{l1}
\end{eqnarray}
The first equality comes from the changes of variables $y =
t/(t+2^k)$ in the first integral and $y = 1/(1+2^kt)$ in the
second integral. As a consequence, the operator $B_2$ can be defined
over the entire Banach space $L^1([0,1])$, and this operator is
continuous with respect to the topology generated by its norm:
$$
B_2 : L^1([0,1]) \longrightarrow L^1([0,1]) \mbox{ and } B_2 \in
\mathcal{L}(L^1([0,1]),L^1([0,1])).
$$
This property was already noticed by Brent in his original article
\cite{brent}. Here is a first application of Lemma \ref{lem}:

\begin{pr}\label{commute}
If $h$ is a function in $K_1$, c.f. (\ref{K1}),
then $F(h)'=B_2(h')$ in $L^1([0,1])$.
\end{pr}

\begin{proof}
Let $h \in S \cap \{ g  |  \, g(0)=1, g(1)= 0 \}$. Let us define
$$
f_n(x) = F_n(h)(x) = \sum_{k = 1}^{n} 2^{-k} \left( h \left(\frac{x}{x+2^k}
\right) - h \left(\frac{1}{1+2^kx} \right) \right).
$$ 
By Lemma \ref{stable} and \ref{conv}, these functions are convex,
decreasing and continuously differentiable over $]0,1]$. Therefore,
over $]0,1]$, we have
$$
f_n'(x) = \sum_{k = 1}^{n} \left(\frac{1}{x+2^k} \right)^2 h'
\left(\frac{x}{x+2^k} \right) + \left(\frac{1}{1+2^kx} \right)^2 h'
\left(\frac{1}{1+2^kx} \right).
$$
The sequence $\{f_n'\}_{n \in \N}$, being the partials sum of the
series that defines $B_2$, converges to $B_2(h')$ in $L^1([0,1])$.
The condition of Lemma \ref{lem} are fulfilled and since the sequence
$\{f_n\}_{n \in \N}$ converges uniformly towards $F(h)$, we have
$F(h)'=B_2(h')$ in $L^1([0,1])$.

In general, if $h \in \overline{S} \cap \{ g | \, g(0)=1, g(1)= 0 \}$,
there exists a sequence $\{h_n\}_{n \in\N}$ of function in $S \cap \{
g | \, g(0)=1, g(1)= 0 \}$ that converges uniformly to $h$. We use
again Lemma \ref{lem} to have that $h_n' \longrightarrow h'$ in
$L^1$. Since every $h_n$ belongs to $S \cap \{ g | \, g(0)=1, g(1)= 0
\}$ the partial result above applies and thus $F(h_n)'=B_2(h_n')$ for
all $n \in \N$. Since $\{F(h_n)\}_{n \in\N}$ is a sequence of
decreasing convex functions that converges uniformly to $F(h)$, thanks
to the continuity of $F$ with respect to $||.||_{\infty}$, we can once
again apply Lemma \ref{lem} to this sequence. Taking the limit leads
to the result since $F(h_n)' \longrightarrow F(h)'$ in $L^1$ and
$B(h_n') \longrightarrow B(h')$ in $L^1$ as well, because of the
continuity of $B_2$ with respect to $||.||_{1}$.
\end{proof}

Based on these properties, we can prove the following expected
theorem, using another time Lemma \ref{lem}:

\begin{thm}\label{thm}
  The sequence $\{g_n\}_{n \in\N}$ possesses a unique accumulation
  point, and therefore converges uniformly to a limiting function
  $g_{\infty}$ which is a fixed point of the linear operator $F$.  The
  sequence of derivatives $\{g'_n\}_{n \in\N}$ converges almost
  everywhere and in the $L^1$-norm to $g'_{\infty}$, which is a fixed
  point of the linear operator $B_2$.
\end{thm}

\begin{proof}
Let us consider $g_{\infty}$, an accumulation point of the sequence
$\{g_n\}_{n \in\N}$. Based on Proposition \ref{commute}, we have
$$
g_{\infty}' = F(g_{\infty})' = B_2(g_{\infty}') \mbox{ in } L^1([0,1]).
$$
Consider the following sequence of non-negative real numbers:
$$
u_n = \int_0^1 |g_{\infty}'-g'_n|dt= || g_{\infty}'-g'_n ||_1 \; ,\; \;
n \in \N.
$$
Then, using the previous equality, Proposition \ref{commute} and
the fact that $B_2$ is a contraction, we see that the sequence
$\{u_n\}_{n \in\N}$ is decreasing:
\begin{eqnarray*}
u_{n+1} & = & ||g_{\infty}'-g'_{n+1}||_1\\
    & = & ||g_{\infty}'-F(g_{n})'||_1\\
    & = & ||B_2(g_{\infty}')-B_2(g_{n}')||_1\\
    & = & ||B_2(g_{\infty}'-g_{n}')||_1\\
    & \leq & ||g_{\infty}'-g_{n}'||_1\\
    & = & u_n.
\end{eqnarray*}
Let $\{g_{n_k}\}_{k \in\N}$ be a subsequence of $\{g_n\}_{n \in\N}$
that converges to $g_{\infty}$. Note that the conditions of Lemma
\ref{lem} are fulfilled and therefore the sequence $\{g'_{n_k}\}_{k
  \in\N}$ converges in $L^1([0,1])$ and almost everywhere to
$g_{\infty}'$.  This implies that the decreasing sequence $\{u_n\}_{n
  \in \N}$ possesses a subsequence that converges to $0$ and therefore
$$
\lim_{n \rightarrow \infty} u_n = 0.
$$
In other words, the sequence $\{g'_{n}\}_{n \in\N}$ converges to
$g_{\infty}'$ in $L^1([0,1])$. Thus, since
$$
\lim_{n \rightarrow \infty} g_n(x) =\lim_{n \rightarrow \infty}
\left( g_n(0) + \int_0^x g_n'(t) dt \right)
= 1+ \int_0^x g_{\infty}'(t) dt = g_{\infty}(x),
$$
we see that the sequence $\{g_n\}_{n \in\N}$ converges
point-wise to $g_{\infty}$. This makes impossible the existence of
another accumulation point. As explained at the beginning of the
section, this shows that the sequence $\{g_n\}_{n \in\N}$
converges to $g_{\infty}$.
\end{proof}

\Section{Conclusion}

Theorem \ref{thm} shows that the operator $B_2$ has a unique eigenfunction
with eigenvalue 1. Computational experiments \cite{brent2} show that the next
eigenvalues seem to be conjugate complex numbers $\lambda$ and
$\overline{\lambda}$ close to $0.1735 \pm 0.00884i$, with
$|\lambda_1|=|\lambda_2|=0.1948$. Therefore, $B_2$ seems to present the
spectral gap Vall\'ee's operator $V_2$ possesses. The method described in this
article does not seem to extend in such a way that this spectral gap can be
proved. This would imply the exponential speed of convergence already
suspected in \cite{brent}.

We did not prove that the function $g_{\infty}'$ is continuous over
$]0,1]$, which is strongly suspected. If proven, this continuity would
directly imply the uniform convergence of the sequence of continuous
and increasing functions $\{g_n'\}_{n \in\N}$ over any compact set of
$]0,1]$ because of a theorem of Dini.

\Sectionn{Acknowledgment}

I would like to thank Richard Brent, Joachim Rosenthal and Brigitte
Vall\'ee for their numerous remarks and advices.


\bibliographystyle{amsplain}

\end{document}